\documentclass{amsart}

\newtheorem{thm}{Theorem}
 \newtheorem{cor}[thm]{Corollary}
 \newtheorem{lem}[thm]{Lemma}
 \newtheorem{prop}[thm]{Proposition}
 \newtheorem{defn}[thm]{Definition}

\newcommand{\ma}{\mathcal{A}}
\newcommand{\mb}{\mathcal{B}}
\newcommand{\mc}{\mathcal{C}}
\newcommand{\md}{\mathcal{D}}

\newcommand{\mm}{\mathcal{M}}

\newcommand{\emc}{E_\mathcal{C}}
\newcommand{\emb}{E_\mathcal{B}}

\newcommand{\kmb}{k_\mathcal{B}}
\newcommand{\kmc}{k_\mathcal{C}}

\begin{document}

\title{Operator-valued free Fisher information of random matrices}

\author{Bin Meng}
\address{Bin Meng\newline College of Science,\newline Nanjing University of
Aeronautics and Astronautics,\newline Nanjing 210016, P.R.China}
\email{b.meng@nuaa.edu.cn}

\thanks{2000 \textit{Mathematics Subject Classification.} Primary 46L54, 42C15}

\date{}

\thanks{\textit{Key words and phrases.} Conjugate variable, free Fisher information,
 random matrix, modular frame.
 }

\begin{abstract}
We study the operator-valued free Fisher information of random
matrices in an operator-valued noncommutative probability space.
We obtain a formula for
$\Phi^\ast_{M_2(\mb)}(A,A^\ast,M_2(\mb),\eta)$, where $A\in
M_2(\mb)$ is a $2\times 2$ operator matrix on $\mb$,  and $\eta$
is linear operators on $M_2(\mb)$. Then we consider a special
setting: $A$ is an operator-valued semicircular matrix with
conditional expectation covariance, and find that
$\Phi_\mb^\ast(c,c^\ast:\mb,id)=2Index(E)$, where $E$ is a
conditional expectation of $\mb$ onto $\md$ and $c$ is a circular
variable with covariance $E$.

\end{abstract}
\maketitle

\section{Introduction and preliminaries}
Free probability theory is a noncommutative probability theory
where the classical concept of independence is replaced by the
notion of "freeness". This theory, due to D.Voiculescu, has very
important applications on operator algebras (see \cite{ge01,voi2})

Originally, a noncommutative probability space is a pair
$(\ma,\tau)$, where $\ma$ is a $C^\ast-$ or von Neumann algebra
and $\tau$ is a state on $\ma$. Free independence is defined in
terms of reduced free product relation given by $\tau$. This
notion was generalized by D.Voiculescu and others to an algebra
valued noncommutative probability space where $\tau$ is replaced
by a conditional expectation $\emb$ onto a subalgebra $\mb$ of
$\ma$, and freeness is replaced by freeness with amalgamation.

\begin{defn}\cite{voi1}
Let $\ma$ be a unital algebra over $\Bbb{C}$, and let $\mb$ be a
subalgebra of $\ma$, $1\in \mb$. $\emb:\, \ma\longrightarrow \mb$
is a conditional expectation, i.e. a linear map such that
$\emb(b_1ab_2)=b_1\emb(a)b_2 ,\,\emb(b)=b,$ for any $b,b_1,b_2 \in
\mb,\,a \in \ma$. We call $(\ma,\emb,\mb)$ an operator-valued (or
$\mb-$ valued) noncommutative probability space and elements in
$\ma$ are called $\mb-$ random variables.

The algebra freely generated by $\mb$ and an indeterminate $X$
will be denoted by $\mb[X]$.
 The distribution of $a\in\ma$ is a conditional expectation
$\mu_a: \,\mb[X] \longrightarrow \mb ,\, \mu_a=\emb \circ \tau_a$,
where $\tau_a:\, \mb[X]\longrightarrow \ma$ is the unique
homomorphism such that $\tau_a(b)=b$, for any $b\in
\mb,\,\tau_a(X)=a$. Let $\mb\subset\ma_i\subset\ma,(i\in I)$
 be subalgebras. The family $(\ma_i)_{i\in I}$ will be called free
 with amalgamation over $\mb$(or $\mb-$free),
 if $\emb(a_1a_2\cdots a_n)=0$ whenever $a_j\in\ma_{i_j}$ with
 $i_1\neq i_2\neq \cdots\neq i_n$ and $\emb(a_j)=0,1\leq j\leq n$.
 A sequence $\{a_i\}_{i\in I}\subseteq\ma$ will be called free with amalgamation over $\mb$ if the family of
 subalgebras generated by $(\mb\bigcup \{a_i\})_{i\in I}$ is
 $\mb-$free.
\end{defn}

In this paper we only consider $\mb-$valued
$W^\ast-$noncommutative probability space, which means $\ma,\mb$
are von Neumann (sub)algebras and $\emb$ is a faithful conditional
expectation (i.e. projection with norm 1) of $\ma$ onto $\mb$.

Classical Fisher information is derived from the statistical
estimation theory which is defined by R.A.Fisher (for some details
we refer to \cite{cov1}). By analogy with the classical case
D.Voiculescu introduced the free Fisher information of
self-adjoint random variables in a tracial $W^\ast-$noncommutative
probability space (see \cite{voi3,voi4}). Then A.Nica,
D.Shlykhteko and R.Speicher investigated it in view of cumulant
function (see \cite{nic1}) and study the free Fisher information
of random matrices (see \cite{nic2}). They found the following
equality holds: $\Phi^\ast(\{a_{ij},a_{ij}^\ast\}_{1\leq i,j\leq
d}:\mb)=d^3\Phi^\ast(\{A,A^\ast\}:M_d(\mb))$, where $\Phi^\ast$
denotes the free Fisher information, $A=[a_{ij}]_{i,j=1}^d\in
M_d(\ma):=\ma\otimes M_d(\Bbb{C})$ is an operator matrix with
entries in $\ma$, and $\mb\subseteq \ma$ is a subalgebra of $\ma$.
D.Shlykhteko even introduced free Fisher information with respect
to a completely positive map (see \cite{shl1}). All of these works
are done in the tracial von Neumann algebra framework. In the case
of states, the non-tracial framework was worked out by D.
Shlyakhtenko in \cite{shl3}. In \cite{meng1,meng2}, we generalized
the notion of free Fisher information to the operator-valued
setting and this work is done in the general von Neumann algebra
framework and found that the free Fisher information of an
operator-valued semicircular variable with conditional expectation
covariance is closely related to modular frames.

In the present paper, based on the notion of operator-valued free
Fisher information with respect to a map which introduced in
\cite{meng2}, we study the free Fisher information of random
matrices and obtain a formula of operator-valued free Fisher
information of random matrices, which generalized Theorem 1.2,
Proposition 4.1 in \cite{nic2}, and the method we prove it is
quite different from that of Theorem 1.2 in \cite{nic2}, since the
conditional expectations are not "tracial" in general.
Furthermore, we find that the free Fisher information of circular
variables is closely related to modular frames and the index of
conditional expectations.

Our main method to prove the results is R.Speicher's cumulant
function technique (see \cite{spe1}), just as in
\cite{meng1,meng2}. We use $\widehat{\kmb}:=(\kmb^{(n)})_{n\geq
1}$ and $\widehat{\emb}:=(\emb^{(n)})$ to denote the cumulant and
moment functions induced by $\emb$ respectively where
$\emb^{(n)}(a_1\otimes\cdots\otimes a_n):=\emb(a_1\cdots a_n)$ and
$\kmb^{(n)}$ is determined by the following recurrence formula
(where $E_\mb^{(0)}=1$ formally):
\begin{eqnarray*}
 \lefteqn{\emb^{(n)}(a_1\otimes\cdots\otimes a_n)=}\\& & =\sum\limits_{r=0}^{n-1}\sum\limits_{1<i(1)<i(2)\cdots<i(r)\leq
n} k^{(r+1)}_\mb(a_1\emb^{(i(1)-2)}(a_1\otimes\\ & &\quad\cdots
\otimes a_{i(1)-1}) \otimes
a_{i(1)}\emb^{(i(2)-i(1)-1)}(a_{i(1)+1}\otimes\cdots\otimes
a_{i(2)-1})\\
& & \quad\otimes a_{i(2)}\cdots\otimes
a_{i(r)}\emb^{(n-i(r))}(a_{i(r)+1}\otimes \cdots\otimes
a_n))\end{eqnarray*}
 for all $n\in\Bbb{N}$ and $a_1,\cdots,a_n\in \ma$. It is easy to see
  $\widehat{\emb}$ and $\widehat{\kmb}$ determine
each other uniquely by the above recurrence formula.

Another tools we use is the frame theory. The frame theory in
Hilbert space is derived from the wavelets theory. D.Larson and
D.Han etc. studied frames by the operator theory (see \cite{han}).
Furthermore, M.Frank and D.Larson generalized the frame notion to
the Hilbert $C^\ast-$module setting called modular frames (see
\cite{fra2,fra3}). For an unital $C^\ast-$algebra $\md$, a Hilbert
$\md-$module is a linear space and algebraic $\md-$module $\mm$
together with a $\md-$valued inner product $\langle.,.\rangle_\md$
and complete with respect to the Hilbert $C^\ast-$module norm:
$\|\langle.,.\rangle_\md\|^{\frac{1}{2}}$. We refer to \cite{lan1}
for more details on Hilbert $C^\ast-$module theory. A sequence
$\{f_j:j\in\Bbb{J}\}\subseteq \mm$ is said to be a frame if there
are real constants $C,D>0$ such that
 $$C\langle x,x\rangle_\md\leq \sum_i\langle x,f_i\rangle_\md
\langle f_i,x\rangle_\md\leq D\langle x,x\rangle_\md, \hbox{ for
all } x\in\mm.$$ If one can choose $C=D=1$ the frame will be
called normalized tight. Note that for a frame the sum in the
middle only converges weakly, in general. If the sum converges in
norm then the frame will be called standard (see \cite{fra2}).
M.~Frank and D.~Larson also proved that every algebraically
finitely or countably generated Hilbert $C^\ast-$module possesses
a standard normalized tight frame (\cite{fra2}). In this paper,
when we say "frames", it means "standard frames".

The following lemma will be used frequently.
\begin{lem}\cite{meng1}
Let $(\ma,\emb,\mb)$ be a noncommutative probability space,
$E:\mb\rightarrow \md$ be a conditional expectation, and let $X\in
(\ma,\emb,\mb)$ be a semicircular variable with covariance $E$.
Denote the Hilbert $C^\ast-$module generated by $\mb$ by
$\mathcal{L}_\md^2(\mb)$ whose inner product induced by $E$ and
let $\mathcal{L}_\md^\ast(\mb)=\mb$. Let $\{f_i\}\subseteq \mb$ be
a frame in $\mathcal{L}_\md^2(\mb)$. Then
\begin{equation}
\Phi^\ast_\mb(X:\mb)=\sum\limits_if_if_i^\ast=Index(E)
\end{equation}
and $X$ is free from $\ma$ with amalgamation over $\md$.
\end{lem}

In the above lemma, $Index(E)$ denotes the index of the
conditional expectation $E$. We refer to \cite{bai1,fra1,wat1} for
more details on the index of conditional expectations.

The paper is organized as follows: In section 2, we review the
definition and some results of the operator-valued free Fisher
information with respect to a linear map. We get a explicit
formula for the operator-valued free Fisher information of random
matrices with respect to a certain kind of linear maps. In section
3, we consider the semicircular random matrices' free Fisher
information. We can compute the free Fisher information of
circular variables with covariance $E$ by modular frames, and
furthermore, we point out that it is just equal to the double of
 index of $E$.

\section{The free Fisher information of random matrices}
In \cite{meng2}, we introduced the operator-valued free Fisher
information of one variable with respect to a linear map which is
a generalization of D.Shlyakhtenko's notion in \cite{shl1}. Such a
notion can be generalized to several random variables'setting
easily.

\begin{defn}
Let $(\ma,\emb,\mb)$ be a $\mb-$valued $W^\ast-$noncommutative
probability space. Suppose that $\mb\subseteq\mc\subseteq\ma$ is a
von Neumann subalgebra and $\{X_i\}_{i=1}^n\subseteq \ma$ is
algebraically free from $\mc$ over $\mb$. Denote by
$\mathcal{L}_\mb^2(\mc[X_1,X_2,\cdots,X_n])$ the Hilbert
$\mb-$module generated by $\mc[X_1,X_2,\cdots,X_n]$, whose inner
product is defined by $\langle x,y\rangle:=\emb(x^\ast y)$, for
any $x,y\in \mc[X_1,X_2,\cdots,X_n]$. Let $\{\eta_i\}_{i=1}^n$ be
a sequence of linear maps on $\mc$. $\{\xi_i\}_{i=1}^n\subseteq
\mathcal{L}_\mb^2(\mc[X_1,X_2,\cdots,X_n])$ will be called the
conjugate system of $\{X_i\}_{i=1}^n$ with respect to
$\mc,\{\eta_i\}_{i=1}^n$, if it satisfying:
\begin{eqnarray}
\lefteqn{\emb(\xi_ic_0X_{i_1}c_1X_{i_2}\cdots X_{i_m}c_m)}\\
&
&=\sum\limits_{j=1}^m\delta_{i,i_j}\emb(\eta_i\emc(c_0X_{i_1}\cdots
X_{i_{j-1}}c_{j-1}))\cdot\emb(c_jX_{i_{j+1}}\cdots
X_{i_m}c_m)\nonumber
\end{eqnarray}
for any $m\geq 0$, $c_0,\cdots,c_m\in\mc$ and
$i_1,\cdots,i_m\in\{1,2,\cdots,m\}$.

If such a $\{\xi_i\}_{i=1}^n$ exists, then we define
$\Phi^\ast_\mb(X_1,\cdots,X_n:\mc,\eta_1,\cdots,\eta_n)$--the free
Fisher information of $\{X_i\}_{i=1}^n$ with respect to
$\mc,\{\eta_i\}_{i=1}^n$ by $\sum\limits_i \xi_i\xi_i^\ast$, that
is,
$$\Phi^\ast_\mb(X_1,\cdots,X_n:\mc,\eta_1,\cdots,\eta_n)=\sum\limits_{i=1}^n
\xi_i\xi_i^\ast$$
\end{defn}

From \cite{meng2,nic1}, we know that the conjugate system can be
expressed by cumulant function.

\begin{prop}
With the above notations, let $\emc$ be a conditional expectation
of $\ma$ onto $\mc$ such that $\emb=\emb\emc$ and let
$(k_\mc^{(n)})_{n\geq 1}$ be the cumulant function induced by
$\emc$. Then $\{\xi_i\}_{i=1}^n\subseteq
\mathcal{L}_\mb^2(\mc[X_1,X_2,\cdots,X_n])$ is the conjugate
system of $\{X_i\}_{i=1}^n$ with respect to
$\mc,\{\eta_i\}_{i=1}^n$ if and only if the following equations
hold
$$\begin{cases}\kmc^{(1)}(\xi_ic)=0, &\forall
i\in\{1,2,\cdots,n\},\\
\kmc^{(2)}(\xi_i\otimes ca)=\delta_{a,X_i}\emb(\eta_i(c)),&\forall
i\in\{1,2,\cdots,n\},\\
\kmc^{(m+1)}(\xi_i\otimes c_1a_1\otimes\cdots\otimes
c_ma_m)=0,&\forall i\in\{1,2,\cdots,n\},m\geq 2\end{cases}$$ where
$c,c_1,\cdots,c_m\in\mc$,
$a,a_1,\cdots,a_m\in\{X_1,\cdots,X_n\}\cup\mc$.
\end{prop}

We introduce some notations about random matrices.
$M_d(\ma):=\ma\otimes M_d(\Bbb{C}), M_d(\mb):=\mb\otimes
M_d(\Bbb{C})$ are the sets of all the $d\times d$ matrices on
$\ma,\mb$ respectively. If $\emb$ is a conditional expectation of
$\ma$ onto $\mb$, then define the conditional expectation
$\emb\otimes I_d: \ma\otimes M_d(\Bbb{C})\rightarrow\mb\otimes
M_d(\Bbb{C})$ by $(\emb\otimes
I_d)([a_{ij}]_{i,j=1}^d):=[\emb(a_{ij})]_{i,j=1}^d$, for all
$A:=[a_{ij}]_{i,j=1}^d \in M_d(\ma)$.  Obviously,
$(M_d(\ma),\emb\otimes I_d,M_d(\mb))$ is a $M_d(\mb)-$valued
noncommutative probability and the elements in it are called
operator-valued random matrices in general.

Since there is no essential difference between $2\times 2$ and
$n\times n$ matrices, for convenience, we only consider $2\times
2$ matrices below.

The following theorem describe the free Fisher information of
random matrices in terms of its entries' free Fisher information.

\begin{thm}
Let $(\ma,\emb,\mb)$ be a $\mb-$valued $W^\ast-$noncommutative
probability space and let $\mb\subseteq \mc\subseteq\ma$ be a
subalgebra. $\{\eta_{ij},\xi_{ij}\}_{i,j=1}^2$ is a sequence of
linear map on $\mc$. Define $\eta: M_2(\mc)\rightarrow M_2(\mc)$
by $\eta\left(\begin{array}{cccc}c_{11} & c_{12}\\ c_{21} &
c_{22}\end{array}\right)=\left(\begin{array}{cccc}\eta_{11}(c_{11})+\eta_{21}(c_{22})
& 0\\ 0 & \eta_{12}(c_{11})+\eta_{22}(c_{22})\end{array}\right)$
and $\xi:M_2(\mc)\rightarrow M_2(\mc)$ by
$\xi\left(\begin{array}{cccc} c_{11} & c_{12}\\ c_{21} &
c_{22}\end{array}\right)=\left(\begin{array}{cccc}\xi_{11}(c_{11})+\xi_{12}(c_{22})
& 0\\ 0 & \xi_{21}(c_{11})+\xi_{22}(c_{22})\end{array}\right)$.
$A:=[a_{ij}]_{i,j=1}^2\in M_2(\ma)$. Let
$\{x_{ij},y_{ij}\}_{i,j=1}^2$ be the conjugate system of
$\{a_{ij},a_{ij}^\ast\}_{i,j=1}^2$ with respect to
$\{\eta_{ij},\xi_{ij}\}_{i,j=1}^2$. Then
$$\Phi^\ast_{M_2(\mb)}(A,A^\ast:M_2(\mc),\eta,\xi)=\left(\begin{array}{cccc}A_{11}
& A_{12}\\ A_{21} & A_{22}\end{array}\right)$$ where
\begin{eqnarray*}
&&A_{11}=\Phi_\mb^\ast(a_{11},a_{11}^\ast:\mc,\eta_{11},\xi_{11})+\emb(x_{21}x_{21}^\ast)+\emb(y_{12}y_{12}^\ast)\\
&&A_{12}=\emb(x_{11}x_{12}^\ast+x_{21}x_{22}^\ast+y_{11}y_{21}^\ast+y_{12}y_{22}^\ast)\\
&&A_{21}=\emb(x_{12}x_{11}^\ast+x_{22}x_{21}^\ast+y_{21}y_{11}^\ast+y_{22}y_{12}^\ast)\\
&&A_{22}=\Phi_\mb^\ast(a_{22},a_{22}^\ast:\mc,\eta_{22},\xi_{22})+\emb(x_{12}x_{12}^\ast)+\emb(y_{21}y_{21}^\ast)
\end{eqnarray*}
\end{thm}

\begin{proof}
Let $X_1:=\left(\begin{array}{cccc}x_{11} & x_{21}\\
x_{12} & x_{22}\end{array}\right)$,
$X_2:=\left(\begin{array}{cccc} y_{11} & y_{12}\\ y_{21} &
y_{22}\end{array}\right)$. Below we show $\{X_1,X_2\}$ is the
conjugate system of $\{A,A^\ast\}$ with respect to $\{\eta,\xi\}$
using Proposition 4.

Denote by $(K^{(n)})_{n\geq 1}$ the cumulant function induced by
$(\emc\otimes I_2)$. Then for any $C_k:=[c_{ij}^{(k)}]_{i,j=1}^2$,
we get the following equalities.
\begin{eqnarray*}
\lefteqn{K^{(1)}(X_1C_1)=(\emc\otimes
I_2)\left(\begin{array}{cccc}x_{11} & x_{21}\\
x_{12} & x_{22}\end{array}\right) \left(\begin{array}{cccc}
c_{11}^{(1)} & c_{12}^{(1)}\\
c_{21}^{(1)} & c_{22}^{(1)}\end{array}\right)}\\
& &=(\emc\otimes
I_2)\left(\begin{array}{cccc}x_{11}c_{11}^{(1)}+x_{21}c_{21}^{(1)}
& x_{11}c_{12}^{(1)}+x_{21}c_{22}^{(1)}\\
x_{12}c_{11}^{(1)}+x_{22}c_{21}^{(1)} &
x_{12}c_{12}^{(1)}+x_{22}c_{22}^{(1)}\end{array}\right)\\
& &=0;
\end{eqnarray*}
By the same way, we can prove $K^{(1)}(X_2C_1)=0$.
\begin{eqnarray*}
\lefteqn{K^{(2)}(X_1\otimes C_1A)}\\
& &=(\emc\otimes I_2)^{(2)}(X_1\otimes C_1A)-(\emc\otimes
I_2)^{(1)}(X_1(\emc\otimes I_2)(C_1A))\\
& &=(\emc\otimes I_2)(X_1C_1A)\\
& &=(\emc\otimes
I_2)\left(\begin{array}{cccc}\sum\limits_{i_1,i_2=1}^2x_{i_11}c_{i_1,i_2}^{(1)}
a_{i_21} &
\sum\limits_{i_1,i_2=1}^2x_{i_11}c_{i_1i_2}^{(1)}a_{i_22}\\
\sum\limits_{i_1,i_2=1}^2x_{i_12}c_{i_1i_2}^{(1)}a_{i_21} &
\sum\limits_{i_1,i_2=1}^2x_{i_12}c_{i_1i_2}^{(1)}a_{i_22}
\end{array}\right)\\
&
&=\left(\begin{array}{cccc}\emb(\eta_{11}(c_{11}^{(1)})+\eta_{21}(c_{22}^{(1)}))
& 0\\
0 &
\emb(\eta_{12}(c_{11}^{(1)})+\eta_{22}(c_{22}^{(1)}))\end{array}\right)\\
& &=(\emb\otimes I_2)(\eta(C_1))
\end{eqnarray*}
Similarly, we get $K^{(2)}(X_2\otimes C_1A^\ast)=\emb(\xi(C_1))$;
$K^{(2)}(X_1\otimes C_1A^\ast)=0$; $K^{(2)}(X_2\otimes C_1A)=0$.

To prove
 \begin{equation}K^{(m+1)}(X_1\otimes C_1A_1\otimes
\cdots\otimes C_mA_m)=0,\forall m\geq 2,\end{equation} where
$A_i\in\{A,A^\ast,M_2(\mc)\}$
 we need to induce on
$m$.

When $m=2$,
\begin{eqnarray*}
\lefteqn{K^{(3)}(X_1\otimes C_1A_1\otimes C_2A_2)}\\
& &=(\emc\otimes I_2)(X_1\cdot C_1A_1\cdot C_2A_2)
   -K^{(2)}(X_1\otimes C_1A_1(\emc\otimes I_2)(C_2A_2))\\
& &\quad -K^{(2)}(X_1\otimes (\emc\otimes I_2)(C_1A_1)C_2A_2)\\
& &=(\emc\otimes I_2)(X_1\cdot C_1A_1\cdot
     C_2A_2)-\delta_{A_1,A}(\emc\otimes I_2)(\eta(C_1))(\emc\otimes
     I_2)(C_2A_2)\\
& &\quad -\delta_{A_2A}(\emc\otimes I_2)\eta((\emc\otimes
      I_2)(C_1A_1)C_2)
\end{eqnarray*}

On the other hand,
\begin{eqnarray*}
\lefteqn{(\emc\otimes I_2)(X_1\cdot C_1A_1\cdot C_2A_2)}\\
& &=(\emc\otimes I_2)\left[\left(\begin{array}{cccc} x_{11} & x_{21}\\
x_{12} &
x_{22}\end{array}\right)\left(\begin{array}{cccc}c_{11}^{(1)} &
c_{12}^{(1)}\\ c_{21}^{(1)} & c_{22}^{(1)}\end{array}\right)
\left(\begin{array}{cccc}a_{11}^{(1)} & a_{12}^{(1)}\\
a_{21}^{(1)} &
a_{22}^{(1)}\end{array}\right)\left(\begin{array}{cccc}
c_{11}^{(2)} & c_{12}^{(2)}\\ c_{21}^{(1)} &
c_{22}^{(2)}\end{array}\right)
 \left(\begin{array}{cccc} a_{11}^{(2)} & a_{12}^{(2)}\\
a_{21}^{(2)} & a_{22}^{(2)}\end{array}\right)\right]\\
&&=\left[\sum\limits_{i_3,i_4,i_5,i_6=1}^2\emc^{(3)}(x_{i_3i_1}\otimes
   (c_{i_3i_4}^{(1)}a_{i_4i_5}^{(1)})\otimes
   (c_{i_5i_6}^{(2)}a_{i_6i_2}^{(2)}))\right]_{i_1,i_2=1}^2\\
&&=\Bigg[\sum\limits_{i_3,i_4,i_5,i_6=1}^2k^{(1)}(x_{i_3i_1}\emc(c_{i_3i_4}^{(1)}a_{i_4i_5}^{(1)}
   c_{i_5i_6}^{(2)}a_{i_6i_2}^{(2)}))\\
&&\quad +k^{(2)}(x_{i_3i_1}\otimes
  c_{i_3i_4}^{(1)}a_{i_4i_5}^{(1)}\emc(c_{i_5i_6}^{(2)}a_{i_6i_2}^{(2)}))\\
&&\quad\quad+k^{(2)}(x_{i_3i_1}\emc(c_{i_3i_4}^{(1)}a_{i_4i_5}^{(1)})\otimes
  c_{i_5i_6}^{(2)}a_{i_6i_2}^{(2)})  +k^{(3)}(x_{i_3i_1}\otimes
c_{i_2i_3}^{(1)}a_{i_3i_4}^{(1)}\otimes
c_{i_4i_5}^{(2)}a_{i_5i_1}^{(2)})\Bigg]_{i_1,i_2=1}^2\\
&&=\delta_{A_1A}(\emb\otimes I_2)(\eta(C_1))(\emc\otimes
I_2)(C_2A_2)+\delta_{A_2A}(\emb\otimes I_2)\eta((\emc\otimes
I_2)(C_1A_1)C_2),
\end{eqnarray*}
so, $K^{(3)}(X_1\otimes C_1A_1\otimes C_2A_2)=0$

Now assume that the Eq.(3) holds for the numbers $\leq m-1$, then
\begin{eqnarray*}
&&\quad K^{(m+1)}(X_1\otimes C_1A_1\otimes\cdots\otimes C_mA_m)\\
&&=(\emc\otimes I_2)(X_1\cdot C_1A_1\cdots
C_mA_m)\\
&&\quad-\sum\limits_{j=1}^{m+1}K^{(2)}(X_1(\emc\otimes
I_2)^{(j-1)}(C_1A_1\otimes\cdots\otimes C_{j-1}A_{j-1})\\
&&\quad \otimes C_jA_j(\emc\otimes
I_2)^{(m-j)}(C_{j+1}A_{j+1}\otimes\cdots\otimes C_mA_m))\\
&&=(\emc\otimes I_2)(X_1\cdot C_1A_1\cdots C_mA_m)\\
&&\quad -\sum\limits_{j=1}^{m+1}\delta_{AA_j}(\emb\otimes
I_2)(\eta(\emc\otimes I_2)(C_1A_1\cdots
C_{j-1}A_{j-1})C_j)\\
&&\quad\quad(\emc\otimes
I_2)^{(m-j)}(C_{j+1}A_{j+1}\otimes\cdots\otimes C_mA_m))\\
&&=\Bigg[\sum\limits_{i_3,i_4,\cdots,i_{2(m+1)}=1}^2\emc(x_{i_3i_1}(c_{i_3i_4}^{(1)}
a_{i_4i_5}^{(1)})\cdots
(c_{i_{2(m+1)-1}}^{(m)}a_{i_{2(m+1)}i_2}^{(m)}))\Bigg]_{i_1,i_2=1}^2\\
&&\quad-\sum\limits_{j=1}^{m+1}\delta_{AA_j}(\emb\otimes
I_2)\eta((\emc\otimes I_2)(C_1A_1\cdots
C_{j-1}A_{j-1})C_j)\\
&&\quad\quad(\emc\otimes
I_2)^{(m-j)}(C_{j+1}A_{j+1}\otimes\cdots\otimes C_mA_m))\\
&&=\Bigg[\sum\limits_{i_3,i_4,\cdots,i_{2(m+1)}}^2\sum\limits_{j=0}^mk^{(2)}
(x_{i_3i_1}\emc(c_{i_3i_4}a_{i_4i_5}^{(1)}\cdots
c^{(j)}_{i_{2(j+1)-1}i_{2(j+1)}}a^{(j)}_{i_{2(j+1)}i_{2(j+1)+1}})\\
&&\quad\otimes
c^{(j+1)}_{i_{2(j+2)-1}i_{2(j+2)}}a^{(j+1)}_{i_{2(j+2)}i_{2(j+2)}}\emc(c^{(j+2)}_{i_{2(j+2)+1}
i_{2(j+2)+2}}a_{i_{2(j+2)+2}i_{2(j+2)+3}}\cdots
c^{(m)}_{i_{2m+1}i_{2m+2}}a^{(m)}_{i_{2m+2}i_2}))\Bigg]_{i_1,i_2=1}^2\\
&&\quad-\sum\limits_{j=1}^{m+1}\delta_{AA_j}(\emb\otimes
I_2)\eta((\emc\otimes I_2)(C_1A_1\cdots
C_{j-1}A_{j-1})C_j)\\
&&\quad\quad(\emc\otimes
I_2)^{(m-j)}(C_{j+1}A_{j+1}\otimes\cdots\otimes C_mA_m))\\
&&=\Bigg[\sum\limits_{i_3,\cdots,i_{2(m+1)}=1}^2\delta_{a_{i_3i_1},a^{(j+1)}_{i_{2(j+2)}i_{2(j+2)}}}
\emb\eta(\emc(c_{i_3i_4}a_{i_4i_5}^{(1)}\cdots
c^{(j)}_{i_{2(j+1)-1}i_{2(j+1)}}a^{(j)}_{i_{2(j+1)}i_{2(j+1)+1}})\\
&&\quad\quad c^{(j+1)}_{i_{2(j+2)-1}i_{2(j+2)}})
\emc(c^{(j+2)}_{i_{2(j+2)+1}
i_{2(j+2)+2}}a_{i_{2(j+2)+2}i_{2(j+2)+3}}\cdots
c^{(m)}_{i_{2m+1}i_{2m+2}}a^{(m)}_{i_{2m+2}i_2}))\Bigg]_{i_1,i_2=1}^2\\
&&\quad-\sum\limits_{j=1}^{m+1}\delta_{AA_j}(\emb\otimes
I_2)\eta((\emc\otimes I_2)(C_1A_1\cdots
C_{j-1}A_{j-1})C_j)\\
&&\quad\quad(\emc\otimes
I_2)^{(m-j)}(C_{j+1}A_{j+1}\otimes\cdots\otimes C_mA_m))\\
&&=0
\end{eqnarray*}
Thus we claim that $X_1,X_2$ is the conjugate system of
$(A,A^\ast)$ with respect to $(\eta,\xi)$. So,
\begin{eqnarray*}
&&\Phi^\ast_\mb(A,A^\ast:M_2(\mc),\eta,\xi)\\
&=&(\emb\otimes I_2)\left[\left(\begin{array}{cccc} x_{11} &
x_{21}\\ x_{12} &
x_{22}\end{array}\right)\left(\begin{array}{cccc}x_{11}^\ast &
x_{12}^\ast\\ x_{21}^\ast &
x_{22}^\ast\end{array}\right)+\left(\begin{array}{cccc} y_{11} &
y_{12}\\ y_{21} &
y_{22}\end{array}\right)\left(\begin{array}{cccc}y_{11}^\ast &
y_{21}^\ast\\ y_{12}^\ast & y_{22}^\ast\end{array}\right)\right]\\
&=&\left(\begin{array}{cccc}A_{11} & A_{12}\\ A_{21} &
A_{22}\end{array}\right)
\end{eqnarray*}
\end{proof}

The following corollary is a analogue of Theorem 1.2 in
\cite{nic2}.
\begin{cor}
Let $\mb_T^{(2)}:=\left\{\left(\begin{array}{cccc}b & 0\\ 0 &
b\end{array}\right)|b\in\mb\right\}$ and $(\emb\otimes
tr_T)\left(\begin{array}{cccc}b_{11} & b_{12} \\ b_{21} &
b_{22}\end{array}\right):=\left(\begin{array}{cccc}\emb\left(\dfrac{b_{11}+b_{22}}{2}\right)
& 0\\ 0 &
\emb\left(\dfrac{b_{11}+b_{22}}{2}\right)\end{array}\right)$. Then
 in the $\mb_T^{(2)}-$valued noncommutative probability space
$(M_2(\ma),\emb\otimes tr_T,\mb_T^{(2)})$, for a random matrix
$A=[a_{ij}]_{i,j=1}^2$ we have
\begin{equation*}
\Phi^\ast_{\mb_T^{(2)}}(A,A^\ast:M_2(\mc),id)=\frac{1}{8}\left(\begin{array}{cccc}\sum\limits_{i,j=1}^2\Phi^\ast_\mb(a_{ij},a_{ij}^\ast:\mc,id)
& 0\\ 0 &
\sum\limits_{i,j=1}^2\Phi^\ast_\mb(a_{ij},a_{ij}^\ast:\mc,id)\end{array}\right)
\end{equation*}
\end{cor}
\begin{proof}
With the notations as in Theorem and its proof, let
$\eta_{ij}=\xi_{ij}=1, \, i,j=1,2$. Then from the proof of
Theorem, we know
\begin{eqnarray*}
K^{(2)}(\dfrac{1}{2} X_1\otimes
CA)&=&\left(\begin{array}{cccc}\dfrac{\emb(c_{11}+c_{22})}{2} &
\quad\\ \quad & \dfrac{\emb(c_{11}+c_{22})}{2}\end{array}\right)\\
&=&(\emb\otimes tr_T)(C)
\end{eqnarray*}
and
\begin{eqnarray*}
K^{(2)}(\dfrac{1}{2} X_2\otimes
CA^\ast)&=&\left(\begin{array}{cccc}\dfrac{\emb(c_{11}+c_{22})}{2}
&
\quad\\ \quad & \dfrac{\emb(c_{11}+c_{22})}{2}\end{array}\right)\\
&=&(\emb\otimes tr_T)(C).
\end{eqnarray*}
Hence $\{\frac{1}{2}X_1,\frac{1}{2}X_2\}$ is the conjugate system
of $A,A^\ast$ with respect to $\mb_T^{(2)},id$ and
\begin{eqnarray*}
&&\Phi^\ast_{\mb_T^{(2)}}(A,A^\ast:M_2(\mc),id)\\
&=&\frac{1}{4}\emb\otimes tr_T(X_1X_1^\ast)+\frac{1}{4}\emb\otimes
tr_T(X_2X_2^\ast)\\
&=&\frac{1}{8}\left(\begin{array}{cccc}\sum\limits_{i,j=1}^2\Phi^\ast_\mb(a_{ij},a_{ij}^\ast:\mc,id)
& 0\\ 0 &
\sum\limits_{i,j=1}^2\Phi^\ast_\mb(a_{ij},a_{ij}^\ast:\mc,id)\end{array}\right)
\end{eqnarray*}

\end{proof}

\section{The free Fisher information of random matrices with
semicircular and circular entries}

Recall that $X\in (\ma,\emb,\mb)$ will be called a semicircular
variable with covariance $\eta$ (or $\eta-$semicircular variable,
where $\eta$ is a linear map on $\mb$), if it satisfies:
$\kmb^{(1)}(X)=0,\, \kmb^{(2)}(X\otimes bX)=\eta(b),\,
\kmb^{(m+1)}(X\otimes b_1X\otimes\cdots\otimes b_mX)=0,\forall
m\geq 2$, where $(\kmb^{(n)})_{n\geq 1}$ is the cumulant function
induced by $\emb$. Let $X,Y$ be two free $\eta-$semicircular
variables, then $C:=\dfrac{X+iY}{\sqrt{2}}$ will be called a
circular variable in $(\ma,\emb,\mb)$ with covariance $\eta$ (or
$\eta-$circular).

In the classical free probability theory, it is well known that
$\left(\begin{array}{cccc}f_1 & c\\ c^\ast &
f_2\end{array}\right)\in (M_2(\ma),\tau\otimes tr, \Bbb{C})$ is
semicircular, where $f_1,f_2,\{c,c^\ast\}$ is a free family in
$(\ma,\tau)$, $f_1,f_2$ are semicircular and $c$ is circular(see
\cite{voi0}). In fact we can get more.
\begin{lem}\cite{shl2}
Let $(\ma,\emb,\mb)$ be a $\mb-$valued noncommutative probability
space, $X,Y\in (\ma,\emb,\mb)$ be $\eta-$semicircular, and let
$C\in (\ma,\emb,\mb)$ be $\eta-$circular. If $X, Y,\{C,C^\ast\}$
is a $\mb-$free family, then $\left(\begin{array}{cccc} X & C\\
C^\ast & Y\end{array}\right)\in (M_2(\ma),\emb\otimes
I_2,M_2(\mb))$ is a semicircular variable with covariance
$\eta^+$, where $\eta^+\left(\begin{array}{cccc} b_{11} & b_{12}\\
b_{21} &
b_{22}\end{array}\right):=\left(\begin{array}{cccc}\eta(b_{11}+b_{22})
& 0\\ 0 & \eta(b_{11}+b_{22})\end{array}\right)$
\end{lem}

\begin{thm}
Let $c\in(\ma,\emb,\mb)$ be $E-$circular, where
$E:\mb\rightarrow\md$ is a conditional expectation. Then
\begin{equation}\Phi^\ast_\mb(c,c^\ast:\mb,E)=2\end{equation}
\end{thm}
\begin{proof}
Let $\{f_i\}\subseteq \mb$ be a normalized tight frame in the
Hilbert $C^\ast-$module $\mathcal{L}_\md^2(\mb)$ which induced by
$E$. Let $s_1,s_2\in (\ma,\emb,\mb)$ be two  $E-$semicircular
variables such that $s_1,s_2,\{c,c^\ast\}$ is a free family with
amalgamation over $\mb$, and let $\{a,x,y,b\}$ be the conjugate
system of $\{s_1,c,c^\ast,s_2\}$ with respect to $E$ Then we claim
$\emb(ax)=\emb(ay)=\emb(xb)=\emb(yb)=0$ since $a,\{x,y\},b$ is a
free family. By Theorem 5,
\begin{eqnarray}
&&\Phi^\ast_{M_2(\mb)}\left(\left(\begin{array}{cccc}s_1 & c\\
c^\ast &
s_2\end{array}\right):M_2(\mb),E^+\right)\\
&=&\left(\begin{array}{cccc}\Phi^\ast_\mb(s_1;\mb,E)+\emb(yy^\ast)
& 0\\ 0& \Phi^\ast_\mb(s_2;\mb,E) &
\emb(xx^\ast)\end{array}\right)\nonumber
\end{eqnarray}
and from \cite{meng2}, we know
$\Phi^\ast_\mb(s_1;\mb,E)=\Phi^\ast_\mb(s_2;\mb,E)=1$.

 On the other hand, $\left(\begin{array}{cccc}s_1 &
c\\ c^\ast & s_2\end{array}\right)\in (M_2(\ma),\emb\otimes
I_2,M_2(\mb))$ is $E^+$semicircular. Note that $E^+$ is not a
conditional expectation, since $E^+\left(\begin{array}{cccc} d &
0\\ 0 & d\end{array}\right)\neq \left(\begin{array}{cccc} d & 0\\
0 & d\end{array}\right)$ in general. But for
 $\Big\{\left(\begin{array}{cccc} f_i & 0\\ 0 &
0\end{array}\right)\Big\}_i\bigcup\Big\{\left(\begin{array}{cccc}
0 & f_i\\ 0 &
0\end{array}\right)\Big\}_i\bigcup\Big\{\left(\begin{array}{cccc} 0 & 0\\
f_i &
0\end{array}\right)\Big\}_i\bigcup\Big\{\left(\begin{array}{cccc}
0 & 0\\ 0 & f_i\end{array}\right)\Big\}_i$, we claim that
$Y:=\sum\limits_i\left(\begin{array}{cccc}f_i & 0\\ 0
&0\end{array}\right)\left(\begin{array}{cccc} s_1 & c\\ c^\ast &
s_2\end{array}\right)E^+\left(\left(\begin{array}{cccc} f_i^\ast &
0\\ 0 &
0\end{array}\right)\right)\\+\sum\limits_i\left(\begin{array}{cccc}0
& 0\\ 0 & f_i\end{array}\right)\left(\begin{array}{cccc} s_1 & c\\
c^\ast & s_2\end{array}\right)E^+\left(\left(\begin{array}{cccc} 0
& 0\\ 0 & f_i^\ast\end{array}\right)\right)$ is the conjugate
variable of $\left(\begin{array}{cccc} s_1 & c\\ c^\ast &
s_2\end{array}\right)$ with respect to $M_2(\mb),E^+$ in
$(M_2(\ma),\emb\otimes I_2,M_2(\mb))$. In fact,
\begin{eqnarray*}
&&K^{(2)}\left(Y\otimes \left(\begin{array}{cccc} b_{11} & b_{12}\\
b_{21} & b_{22}\end{array}\right)\left(\begin{array}{cccc}s_1 &
c\\ c^\ast & s_2\end{array}\right)\right)\\
&=&\sum\limits_i\left(\begin{array}{cccc}f_i & 0\\ 0
&0\end{array}\right)E^+\left(\begin{array}{cccc}f_i^\ast & 0\\ 0 &
0\end{array}\right)E^+\left(\begin{array}{cccc} b_{11} & b_{12}\\
b_{21} & b_{22}\end{array}\right)\\
&&\quad +\sum\limits_i \left(\begin{array}{cccc}0 & f_i\\ 0
&0\end{array}\right)E^+\left(\begin{array}{cccc} 0 & 0\\ f_i^\ast & 0\end{array}\right)
E^+\left(\begin{array}{cccc} b_{11} & b_{12}\\
b_{21} & b_{22}\end{array}\right)\\
&&\quad +\sum\limits_i\left(\begin{array}{cccc} 0 & 0\\ f_i &
0\end{array}\right)E^+\left(\begin{array}{cccc} 0 & f_i^\ast \\ 0
& 0\end{array}\right)E^+\left(\begin{array}{cccc}b_{11} & b_{12}\\
b_{21} & b_{22}\end{array}\right)\\
&&\quad +\sum\limits_i \left(\begin{array}{cccc}0 & 0\\ 0 &
f_i\end{array}\right)E^+\left(\begin{array}{cccc}0 & 0\\ 0 &
f_i^\ast\end{array}\right)E^+\left(\begin{array}{cccc}b_{11} &
b_{12}\\ b_{21} & b_{22}\end{array}\right)\\
&=&\left(\begin{array}{cccc} E(b_{11}+b_{22}) & 0\\ 0 &
E(b_{11}+b_{22})\end{array}\right)=E^+\left(\begin{array}{cccc}b_{11}
& b_{12} \\ b_{21} & b_{22}\end{array}\right),
\end{eqnarray*}
and
$$K^{(m+1)}\left(Y\otimes B_1\left(\begin{array}{cccc}s_1 & c\\ c^\ast
& s_2\end{array}\right)\otimes\cdots\otimes B_m
\left(\begin{array}{cccc} s_1 & c\\ c^\ast &
s_2\end{array}\right)\right)=0,$$ for all $B_1,\cdots, B_m\in
M_2(\mb)$, $m\neq 1$, since $\left(\begin{array}{cccc} s_1 & c\\
c^\ast & c^\ast\end{array}\right)$ is semicircular.

So,
\begin{eqnarray*}
&&\Phi_{M_2(\mb)}^\ast\left(\left(\begin{array}{cccc}s_1 & c\\
c^\ast & s_2\end{array}\right):M_2(\mb),E^+\right)\\
&=&(\emb\otimes I_2)(YY^\ast)\\
&=&(\emb\otimes
I_2)\Bigg[\sum\limits_{i,j}\left(\begin{array}{cccc} f_i & 0\\ 0 &
0\end{array}\right)\left(\begin{array}{cccc}s_1 & c\\ c^\ast &
s_2\end{array}\right)E^+\left(\begin{array}{cccc} f_i^\ast & 0\\ 0
& 0\end{array}\right)E^+\left(\begin{array}{cccc} f_j & 0\\ 0 &
0\end{array}\right)\left(\begin{array}{cccc} s_1 & c\\ c^\ast &
s_2\end{array}\right)\left(\begin{array}{cccc}
f_j^\ast & 0\\ 0 & 0\end{array}\right)\\
&&\quad \sum\limits_{i,j}\left(\begin{array}{cccc} 0 & f_i \\ 0 &
0\end{array}\right)\left(\begin{array}{cccc} s_1 & c\\ c^\ast &
s_2\end{array}\right) E^+\left(\begin{array}{cccc} 0 & 0\\
f_i^\ast & 0\end{array}\right)E^+\left(\begin{array}{cccc} 0 &
f_i\\ 0 &0\end{array}\right)\left(\begin{array}{cccc} s_1 & c\\
c^\ast & s_2\end{array}\right)\left(\begin{array}{cccc} 0 & 0\\
f_i^\ast & 0\end{array}\right)\\
&&\quad +\sum\limits_{i,j}\left(\begin{array}{cccc} 0 & 0\\ f_i &
0\end{array}\right)\left(\begin{array}{cccc}s_1 & c\\ c^\ast &
s_2\end{array}\right)E^+\left(\begin{array}{cccc}0 & f_i^\ast \\ 0
&0\end{array}\right)E^+\left(\begin{array}{cccc} 0 &0\\ f_i &
0\end{array}\right)\left(\begin{array}{cccc} s_1 & c\\ c^\ast &
s_2\end{array}\right)\left(\begin{array}{cccc} 0 & f_i^\ast\\ 0 &
0\end{array}\right)\\
&&\quad+\sum\limits_{i,j}\left(\begin{array}{cccc} 0 & 0\\ 0 &
f_i\end{array}\right)\left(\begin{array}{cccc} s_1 & c\\ c^\ast &
s_2\end{array}\right)E^+\left(\begin{array}{cccc}0 & 0\\ 0 &
f_i^\ast\end{array}\right)E^+\left(\begin{array}{cccc} 0 & 0\\ 0 &
f_i\end{array}\right)\left(\begin{array}{cccc} s_1 & c\\ c^\ast &
s_2\end{array}\right)\left(\begin{array}{cccc} 0 & 0\\ 0 &
f_i^\ast\end{array}\right)\Bigg]\\
&=&(\emb\otimes
I_2)\Bigg[\sum\limits_{i,j}\left(\begin{array}{cccc}f_i & 0\\ 0 &
0\end{array}\right)E^+\left(E^+\left(\begin{array}{cccc}f_i^\ast &
0\\ 0 & 0\end{array}\right)E^+\left(\begin{array}{cccc}f_j & 0\\ 0
& 0\end{array}\right)\right)\left(\begin{array}{cccc} f_j^\ast &
0\\ 0 & 0\end{array}\right)\\
&&\quad \sum\limits_{i,j}\left(\begin{array}{cccc} 0 & f_i\\ 0 &
0\end{array}\right)E^+\left(E^+\left(\begin{array}{cccc} 0 & 0\\
f_i^\ast & 0\end{array}\right)E^+\left(\begin{array}{cccc} 0 & f_i\\
0 & 0\end{array}\right)\right)\left(\begin{array}{cccc}  0 &
0\\ f_i^\ast & 0\end{array}\right)\\
&&\quad+\sum\limits_{i,j}\left(\begin{array}{cccc}0 & 0\\ f_i &
0\end{array}\right)E^+\left(E^+\left(\begin{array}{cccc} 0 &
f_i^\ast\\ 0 & 0\end{array}\right)E^+\left(\begin{array}{cccc} 0 &
\\ f_i & 0\end{array}\right)\right)\left(\begin{array}{cccc} 0 &
f_i^\ast\\ 0 & 0\end{array}\right)\\
&&\quad+\sum\limits_{i,j}\left(\begin{array}{cccc} 0 & 0\\ 0 &
f_i\end{array}\right)E^+\left(E^+\left(\begin{array}{cccc}0 & 0\\
0 & f_i\end{array}\right)\right)\left(\begin{array}{cccc}0 &
f_i^\ast\\ 0 & 0\end{array}\right)\\
&&\quad+\sum\limits_{i,j}\left(\begin{array}{cccc}0 & 0\\ 0 &
f_i\end{array}\right)E^+\left(E^+\left(\begin{array}{cccc} 0 & 0\\
0 & f_i^\ast\end{array}\right)E^+\left(\begin{array}{cccc} 0 & 0\\
0 & f_i\end{array}\right)\right)\left(\begin{array}{cccc} 0 & 0\\
0 & f_i^\ast\end{array}\right)\Bigg]\\
&=&\left(\begin{array}{cccc} 2 & \\ & 2\end{array}\right)
\end{eqnarray*}
Compared with Eq.(6), we get $\emb(yy^\ast)=\emb(xx^\ast)=1$ and
so
$$\Phi_\mb(c,c^\ast:\mb,E)=\emb(yy^\ast)+\emb(xx^\ast)=2$$
\end{proof}

\begin{thm}
Let $c\in(\ma,\emb,\mb)$ be $E-$circular, where
$E:\mb\rightarrow\md$ is a conditional expectation. Then
\begin{equation}\Phi^\ast_\mb(c,c^\ast:\mb,id)=2Index(E)\end{equation}
\end{thm}
\begin{proof}
Let $\{f_i\}\subseteq \mb$ be a normalized tight frame in the
Hilbert $C^\ast-$module $\mathcal{L}_\md^2(\mb)$ which induced by
$E$. Let $s_1,s_2\in (\ma,\emb,\mb)$ be two  $E-$semicircular
variables such that $s_1,s_2,\{c,c^\ast\}$ is a free family with
amalgamation over $\mb$, and let $\{a,x,y,b\}$ be the conjugate
system of $\{s_1,c,c^\ast,s_2\}$ with respect to $E$. Then
\begin{eqnarray}
&&\Phi^\ast_{M_2(\mb)}\left(\left(\begin{array}{cccc}s_1 & c\\
c^\ast &
s_2\end{array}\right):M_2(\mb),id^+\right)\\
&=&\left(\begin{array}{cccc}\Phi^\ast_\mb(s_1;\mb,id)+\emb(yy^\ast)
& 0\\ 0& \Phi^\ast_\mb(s_1;\mb,id) +
\emb(xx^\ast)\end{array}\right)\nonumber
\end{eqnarray}
Then similar to the proof of the above theorem, we can say that\\
$\sum\limits_i\left(\begin{array}{cccc}f_i & o\\ 0 &
0\end{array}\right)\left(\begin{array}{cccc} s_1 & c\\ c^\ast &
s_2\end{array}\right)\left(\begin{array}{cccc} f_i^\ast & \\ &
f_i^\ast\end{array}\right)+\sum\limits_i
\left(\begin{array}{cccc}0 & 0\\ 0 &
f_i\end{array}\right)\left(\begin{array}{cccc} s_1 & c\\ c^\ast &
s_2\end{array}\right)\left(\begin{array}{cccc}f_i^\ast & 0\\ 0 &
f_i^\ast\end{array}\right)$ is the conjugate variable of
$\left(\begin{array}{cccc} s_1 & c\\ c^\ast &
s_2\end{array}\right)$ with respect to $M_2(\mb), id^+$. And so

\begin{eqnarray*}
&&\Phi^\ast_{M_2(\mb)}(\left(\begin{array}{cccc} s_1 & c\\ c^\ast
&
s_2\end{array}\right):M_2(\mb),id^+)=\left(\begin{array}{cccc}2\sum\limits_i
f_if_i^\ast & \\ & 2\sum\limits_i f_if_i^\ast\end{array}\right)\\
&=&\left(\begin{array}{cccc}2Index(E)& \\
&2Index(E)\end{array}\right)
\end{eqnarray*}

Then $\Phi^\ast_\mb(c,c^\ast:\mb,id)=2Index(E)$.

\end{proof}

\bibliographystyle{amsplain}

\end{document}